\documentclass[reqno, 12pt]{amsart}
\usepackage[all]{xy}
 \RequirePackage{amsgen,color}
\RequirePackage{amstext} \RequirePackage{amsbsy}
\RequirePackage{amsopn} \RequirePackage{amsthm}
\RequirePackage{amssymb} \RequirePackage{epsfig}
\RequirePackage{amscd}

\usepackage[mathscr]{eucal}

\newtheorem{lemma}{Lemma}
\newtheorem{cor}{Corollary}

\newtheorem{theorem}{Theorem}

\theoremstyle{definition}

\newtheorem{conj}{Conjecture}[section]

\theoremstyle{remark}

  \let\g\gamma

\def\Tr{{\rm Tr}}
\def\Gal{{\rm Gal}}
\def\Frob{{\rm Frob}}

\def\Z{{\mathbb Z}}

\def\Q{{\mathbb Q}}

\def\G{\Gamma}
\def\ord{{\rm ord}}
\def\prim{{\rm prim}}
\def\SL{SL_2(\mathbb Z)}

\title{Fourier coefficients of noncongruence cuspforms}

\author{Wen-Ching Winnie Li}
\address{Department of Mathematics \\Pennsylvania State University\\
University Park, PA 16802\\USA~and National Center for Theoretical Sciences, Mathematics Division,
National Tsing Hua University, Hsinchu 30013, Taiwan, R.O.C.} \email{wli@math.psu.edu}

\author{Ling Long}
\address{Department of Mathematics\\Iowa State University\\ Ames, IA 50011
\\USA}
\email{linglong@iastate.edu}
\date{}

\thanks{}

\begin{document}

\begin{abstract}
 Given a finite index subgroup of $SL_2(\Z)$ with modular curve defined over $\Q$, under the assumption that the space of weight $k$ ($ \ge 2$) cusp forms is $1$-dimensional, we show
 that a form in this space with Fourier coefficients in $\Q$ has bounded denominators if and only if it is a congruence modular form.
\end{abstract}

\thanks{The research of the first author is partially supported by the NSF grant DMS-0801096, and the second author by the NSF grant DMS-1001332. This project started when the second author visited the first author at the National Center for Theoretical Sciences (NCTS) in Taiwan. She would like to thank NCTS for its hospitality.}
\maketitle
\section{Introduction}

Let $\G$ be a finite index subgroup of $SL_2(\Z)$.  The space of cusp forms of weight $k$ for $\G$, denoted by $S_k(\G)$, consists
of holomorphic functions $f$ on the  Poincar\'e upper half plane $\frak H$
satisfying
\begin{itemize}
\item[(1)] $f| \gamma = f \quad \text{for~all}~~ \gamma \in \G$, { where $|\gamma$ stands for the standard stroke operator  (cf. \cite{shim1})};

\item[(2)] $f$ is holomorphic at all cusps of $\G$;

\item[(3)] $f$ vanishes at all cusps of $\G$.
\end{itemize}

\noindent Suppose the cusp $\infty$ of $\G$ has width $\mu$. A form $f \in S_k(\G)$ has Fourier expansions
of the form $$f(z)=\sum_{n\ge 1} a(n)q^{n/\mu}, ~~~~\text{where}~~~~q = e^{2 \pi iz}.$$
When all Fourier coefficients $a(n)$ are algebraic, we say that $f$ has bounded denominators if there exists a nonzero algebraic integer $c$ such that $c\cdot a(n)$ is algebraically integral for all $n$.

For a congruence subgroup $\G$, the space $S_k(\G)$ has a basis consisting of forms with integral Fourier coefficients (cf. \cite{shim1}), resulting from the fact that Hecke operators have integral eigenvalues. Consequently, any form in $S_k(\G)$ with algebraic coefficients has bounded denominators. When $\G$ is a noncongruence subgroup, which is majority, the situation is unknown. However, a folklore conjecture states that the bounded denominator property should characterize the congruence forms (cf. \cite{birchb94}). More precisely,

\begin{conj}[Unbounded denominators conjecture]\label{conj} A meromorphic weight $k$  modular form  for a finite index subgroup
$\G$ of $SL_2(\Z)$ that is holomorphic on  $\frak H$ (such as the modular $j$ function) with algebraic Fourier coefficients has bounded denominators if and only if it is a cusp form for a congruence subgroup.
\end{conj}

This conjecture has impacts beyond the fundamental developments of modular forms. For instance, it implies a conjecture in the theory of vertex
operator algebras known and believed by physicists, which says that the graded dimension of any
$C_2$-cofinite, rational vertex operator algebra over $\mathbb C$ is a
congruence modular function (cf.~\cite{zhu-96-m-inv,Coste-Gannon99,DLM00}).

\medskip

Given a noncongruence subgroup $\G$, denote by $\G^c$ its congruence closure, i.e., the smallest
congruence subgroup containing $\G$. Suppose $\G^c = \cup_\gamma \G \gamma$. Then $f \mapsto \sum_\gamma f| \gamma$ defines the trace map $$ \Tr^{\G^c}_{\G} : S_k(\G) \to S_k(\G^c).$$   As shown in Berger \cite{berger} (see also \cite{schHecke}), $S_k(\G)$ decomposes into the direct sum of $S_k(\G^c)$ and the kernel $S_k^{\prim}(\G)$ of $\Tr^{\G^c}_{\G}$. The forms in $S_k^{\prim}(\G)$ are thus genuinely noncongruence,   called primitive noncongruence in \cite{schHecke}. The conjecture above then says that any nonzero form in $S_k^{\prim}(\G)$ and hence in $S_k(\G) \smallsetminus S_k(\G^c)$ with algebraic Fourier coefficients must have unbounded denominators.

The purpose of this note is to give a partial supportive answer to the above conjecture.

\begin{theorem}\label{thm:main} { Suppose that the modular curve $X_\G$ of $\G$ has a model defined over $\Q$ so that the cusp at $\infty$ is $\Q$-rational,  $k\ge 2$ and $S_k(\G)$ is 1-dimensional.} Then a form in $S_k(\G)$ with Fourier coefficients in $\Q$ has bounded denominators if and only if it is a congruence modular form.
\end{theorem}

 We exhibit two special occasions to which the above theorem can be applied. The first is when the modular curve $X_\G$ is an elliptic curve defined over $\Q$, the space $S_2(\G)$ is 1-dimensional generated by a form from a holomorphic 1-form on $X_\G$. The second is when an elliptic modular surface fibred over $X_\G$ is a K3 surface defined over $\Q$, $S_3(\G)$ is generated by a form arising from a holomorphic 2-form on the elliptic modular surface.

\medskip

There is another approach to the unbounded denominators conjecture using modular functions and $p$-adic arguments,  see \cite{Kurth-Long06,Kurth-Long08II}.

\section{Galois representations attached to $S_k(\G)$}

Recall that for $k\ge 2$ to the $1$-dimensional space $S_k(\G)$ there are associated compatible family of $2$-dimensional $\ell$-adic representations $\rho_{\G, k, \ell}$ of $\Gal(\overline{\Q}/\Q)$. When $k=2$ this is the dual of the Tate modules of the Jacobian of $X_\G$. For $k \ge 4$ even this is constructed by Scholl in \cite{sch85b}. As remarked in \cite{sch85b} the same should hold for $k \ge 3$ odd, which we shall assume.

 At a prime $p \ne \ell$ where $\rho_{\G, k, \ell}$ is unramified, the characteristic polynomial of $\rho_{\G, k, \ell}$  { of the Frobenius $\Frob_p$ at $p$} has the form $T^2 - b(p)T + \chi(p)p^{k-1}$, where $\chi$ is a real Dirichlet character of $\Z$ and $b(p) \in \Z$  with complex absolute value $|b(p)| \le 2 p^{k-1}$ (cf. \cite{sch85b}).

It follows from the former Serre's conjecture that when $\ell$ is a large prime, the dual of $\rho_{\G, k, \ell}$ is modular (cf. \cite{ALLL}). In other words, there is a normalized newform $g$ of weight $k$ level $N$ and character $\chi$ such that its eigenvalue with respect to the Hecke operator $T_p$ at any prime $p \nmid N$ is $b(p)$. Write $g$ in its Fourier
expansion $g(z) = \sum_{n \ge 1} b(n)q^n$.

Let $f = \sum_{n\ge 1} a(n)q^{n/\mu}$ be a nonzero form in $S_k(\G)$ with rational Fourier coefficients. By Scholl \cite[Theorem 5.4]{sch85b}, $f$ satisfies the congruences predicted by Atkin and Swinnerton-Dyer { (ASD)} \cite{ASD}, namely, there is an integer $N_1$ such that for  every prime $p > N_1$
\begin{equation}\label{eq:asd} a(np)-b(p)a(n)+ \chi(p) p^{k-1} a(n/p)\equiv 0 \mod p^{(k-1)(1+\ord_p n)}, \quad \forall n\ge 1.\end{equation}
When this happens, $f$ and $g$ are said to satisfy the Atkin and Swinnerton-Dyer congruence relations in \cite{lly05}.

We shall prove the theorem below, from which Theorem \ref{thm:main} will follow as $S_k(\G)$ is $1$-dimensional.

\begin{theorem}\label{thm:main2} Suppose $k\ge 2$ and $f(z) = \sum_{n\ge 1} a(n)q^{n/\mu}$ is a form in $S_k(\G)$ with Fourier coefficients in $\Q$ satisfying (1) with a normalized weight $k$ newform $g(z) = \sum_{n \ge 1} b(n)q^n$ for all primes $p > N_1$. If $f$ has bounded denominators, then $f \notin S_k^{\prim}(\G)$.
\end{theorem}

Since $f$ in Theorem \ref{thm:main2} has bounded denominators, we may assume that all Fourier coefficients $a(n)$ are in $\Z$.

\section{Comparing Fourier coefficients of $f$ and $g$}

We begin by observing a bound on Fourier coefficients of $f$.

\begin{theorem}[Selberg \cite{sel65}] Let $f(z)=\sum a(n) q^{n/\mu}$ be a weight $k$ cusp form for some  finite index subgroup of the modular group. Then there exists a constant $C>1$ depending on $f$ such that the complex absolute value of $a(n)$ satisfies the bound
\begin{equation}\label{eq:an}
  |a(n)|<Cn^{k/2-1/5} ~~ ~~\text{for~all} \ ~n\ge 1.
\end{equation}
\end{theorem} Since the coefficients of $g$ satisfy the Ramanujan bound which implies, for any $\varepsilon>0$, $|b(n)|=O(n^{(k-1)/2+\varepsilon})$  for all $n\ge 1$, we may further assume that
\begin{equation}\label{eq:bn}
  |b(n)|<Cn^{k/2-1/5}~~~~ \text{for~all}~~ \ n\ge 1.
\end{equation}

\noindent Suppose $a(m) \ne 0$ for some positive integer $m$. Let $A(m) = m^{k/2-1/5}$. Fix an integer
$P(m)>N_1$ such that $$3C^2 A(m)^2 n^{k/2-1/5}<n^{k-1}$$
 for every integer $n>P(m)$. We proceed to compare $a(mn)$ with $b(n)$. In what follows we use the convention that $a(x) = 0$ if $x$ is not a positive integer.

\begin{lemma}\label{lem:1}
  For primes $p>P(m) $ and all integers $n\ge 1$, we have $a(mp^n)=b(p)a(mp^{n-1})-\chi(p)p^{k-1} a(mp^{n-2})$.  Further, if $p$ is coprime to $m$, we have \begin{equation}\label{eq:4a}
    a(mp^n)=b(p^n)a(m).
  \end{equation}
\end{lemma}
\begin{proof}
  Combining \eqref{eq:an} and \eqref{eq:bn} we have for all $n\ge 0$
  $$|a(mp^{n+1})-b(p)a(mp^n)+\chi(p)p^{k-1}a(mp^{n-1})|< 3C m^{(k/2-1/5)(n+1)} p^{(k/2-1/5)(n+1)} < p^{(k-1)(n+1)}.$$
  The ASD congruence (1) gives  $$a(mp^{n+1})-b(p)a(mp^n)+\chi(p)p^{k-1}a(mp^{n-1})\equiv 0 \mod p^{(k-1)(n+1)}.$$
  It follows from the choice of $P(m)$ that for a prime $p>P(m)$, the integer
  $a(mp^{n+1})-b(p)a(mp^n)+\chi(p)p^{k-1}a(mp^{n-1})$ has to be zero for every $n\ge 0$. Next assume $(p, m) = 1$. When $n=0$, this gives $a(mp)=b(p)a(m)$. We prove the statements by induction. Assume that $a(mp^n)=b(p^n)a(m)$ when $n < e$. Then
 \begin{eqnarray*}
    a(mp^e)&=&b(p)a(mp^{e-1})-\chi(p)p^{k-1}a(mp^{e-2})\\&=&a(m)[b(p)b(p^{e-1})-\chi(p)p^{k-1}b(p^{e-2})]=a(m)b(p^e).
  \end{eqnarray*}
\end{proof}

Next we prove a similar statement with multiple large prime factors.
\begin{lemma}\label{lem:2}Let $p_1,..., p_r$ be $r$ distinct primes greater than $P(m)$ and coprime to $m$. For all integers $e_i \ge 0$ we have
\begin{equation}\label{eq:4}
  a(mp_1^{e_1}\cdots p_r^{e_r})=a(m)b(p_1^{e_1})\cdots b(p_{r-1}^{e_{r-1}})b(p_r^{e^r})= a(m)b(p_1^{e_1}\cdots p_r^{e_r}).
\end{equation}
\end{lemma}
 \begin{proof}  This is Lemma 1 when $r=1$, hence we assume $r \ge 2$.
 We will first use mathematical induction on the sum $e_1+e_2+\cdots +e_r$, and for fixed sum by induction on the number of nonzero $e_i$'s, to show that \begin{equation*}
  a(m)a(mp_1^{e_1}\cdots p_r^{e_r})=a(mp_1^{e_1}\cdots p_{r-1}^{e_{r-1}})a(mp_r^{e_r}).\end{equation*} In view of Lemma \ref{lem:1}, we may assume all $e_i \ge 1$ to begin with. Applying congruence (1) to prime $p_1$, we have
   \begin{equation}\label{eq:5}
   a(mp_1^{e_1}\cdots p_r^{e_r})-b(p_1)a(mp_1^{e_1-1}p_2^{e_2}\cdots p_r^{e_r})+\chi(p_1)p_1^{k-1}a(mp_1^{e_1-2}p_2^{e_2}\cdots p_r^{e_r})\equiv 0 \mod p_1^{(k-1)e_1}
   \end{equation} and
   \begin{multline}\label{eq:6}
     a(mp_1^{e_1}\cdots p_{r-1}^{e_{r-1}})-b(p_1)a(mp_1^{e_1-1}p_2^{e_2}\cdots p_{r-1}^{e_{r-1}})+\chi(p_1)p_1^{k-1}a(mp_1^{e_1-2}p_2^{e_2}\cdots p_{r-1}^{e_{r-1}})\\\equiv 0 \mod p_1^{(k-1)e_1}.
   \end{multline}  Now multiplying \eqref{eq:5} by $a(m)$, we obtain
   \begin{multline*}a(m)a(mp_1^{e_1}\cdots p_{r}^{e_{r}})-a(m)b(p_1)a(mp_1^{e_1-1}p_2^{e_2}\cdots p_{r}^{e_{r}})+\chi(p_1)p_1^{k-1}a(m)a(mp_1^{e_1-2}p_2^{e_2}\cdots p_{r}^{e_{r}})\\ \equiv 0 \mod p_1^{(k-1)e_1}.\end{multline*}  By induction assumption this can be rewritten as
   \begin{multline}\label{eq:7}
a(m)a(mp_1^{e_1}\cdots p_{r}^{e_{r}})-b(p_1)a(mp_1^{e_1-1}p_2^{e_2}\cdots p_{r-1}^{e_{r-1}})a( mp_{r}^{e_{r}})+\\
   \chi(p_1)p_1^{k-1}a(mp_1^{e_1-2}p_2^{e_2}\cdots p_{r-1}^{e_{r-1}}) a(mp_{r}^{e_{r}}) \equiv 0 \mod p_1^{(k-1)e_1}.   \end{multline}
   So  multiplying  both sides of \eqref{eq:6}  by $a(mp_r^{e_r})$ and subtracting it from \eqref{eq:7}, we get
   $$a(m)a(mp_1^{e_1}\cdots p_r^{e_r})-a(mp_1^{e_1}\cdots p_{r-1}^{e_{r-1}})a(mp_r^{e_r})\equiv 0 \mod (p_1^{e_1})^{(k-1)}.$$
   By symmetry, the same congruence holds when one replaces $p_1$ by $p_i$  and $e_1$ by $e_i$ for $1\le i\le r-1$. As $p_1,..., p_{r-1}$ are distinct primes, this gives
   \begin{eqnarray}\label{eq:9}
   a(m)a(mp_1^{e_1}\cdots p_r^{e_r})-a(mp_1^{e_1}\cdots p_{r-1}^{e_{r-1}})a(mp_r^{e_r})\equiv 0  \mod (p_1^{e_1}\cdots p_{r-1}^{e_{r-1}})^{k-1}.
   \end{eqnarray}

   On the other hand, similar to \eqref{eq:5}, we have
   \begin{equation*}
   a(mp_1^{e_1}\cdots p_r^{e_r})-b(p_r)a(mp_1^{e_1}\cdots p_{r-1}^{e_{r-1}}p_r^{e_r-1})+\chi(p_r)p_r^{k-1}a(mp_1^{e_1}\cdots p_{r-1}^{e_{r-1}} p_r^{e_r-2})\equiv 0 \mod p_r^{(k-1)e_r}.
   \end{equation*}
   Multiply $a(m)$ to the above and use induction assumption again to get
   \begin{multline*}
     a(m) a(mp_1^{e_1}\cdots p_r^{e_r})-b(p_r)a(mp_1^{e_1}\cdots p_{r-1}^{e_{r-1}})a(mp_r^{e_r-1})\\+\chi(p_r)p_r^{k-1}a(mp_1^{e_1}\cdots p_{r-1}^{e_{r-1}})a(m p_r^{e_r-2}) \equiv 0 \mod p_r^{(k-1)e_r}.
 \end{multline*}
   Apply $a(mp_r^{e^r})=b(p_r)a(mp_r^{e^{r-1}})-\chi(p_r)p_r^{k-1}a(mp_r^{e^{r-2}})$ from Lemma \ref{lem:1} to the above congruence to obtain
   $$a(m)a(mp_1^{e_1}\cdots p_r^{e_r})-a(mp_1^{e_1}\cdots p_{r-1}^{e_{r-1}})a(mp_r^{e_r})\equiv 0 \mod (p_r^{e_r})^{k-1}.$$
   Combined with \eqref{eq:9}, this gives $$a(m)a(mp_1^{e_1}\cdots p_r^{e_r})-a(mp_1^{e_1}\cdots p_{r-1}^{e_{r-1}})a(mp_r^{e_r})\equiv 0  \mod (p_1^{e_1}\cdots p_{r}^{e_{r}})^{k-1}.$$
   By the choice of $C$, \begin{multline*}
     | a(m)a(mp_1^{e_1}\cdots p_r^{e_r})-a(mp_1^{e_1}\cdots p_{r-1}^{e_{r-1}})a(mp_r^{e_r})|<2C^2A(m)^2(p_1^{e_1}\cdots p_r^{e_r})^{k/2-1/5} \\< (p_1^{e_1}\cdots p_{r}^{e_{r}})^{k-1}.
   \end{multline*}Since $a(m)a(mp_1^{e_1}\cdots p_r^{e_r})-a(mp_1^{e_1}\cdots p_{r-1}^{e_{r-1}})a(mp_r^{e_r})$ is an integer, we conclude that
   $$a(m)a(mp_1^{e_1}\cdots p_r^{e_r})=a(mp_1^{e_1}\cdots p_{r-1}^{e_{r-1}})a(mp_r^{e_r}).$$ By induction hypothesis, this implies\begin{equation*}
  a(m)^{r-1}a(mp_1^{e_1}\cdots p_r^{e_r})=a(mp_1^{e_1})\cdots a(mp_{r-1}^{e_{r-1}})a(mp_r^{e^r}).
\end{equation*} This is the same as \eqref{eq:4} because of Lemma \ref{lem:1}.
\end{proof}

Consequently, $a(mM) = a(m)b(M)$ as long as all prime factors of $M$ are larger than $P(m)$ and coprime to $m$.

\section{Proof of Theorem \ref{thm:main2}}

Recall that a form $h \in S_k(\G)$ is primitive noncongruence if $\Tr_\G^{\G^c} h = 0$. We study the primitive noncongruence condition. The lemma below says that the primitive noncongruence property is independent of the choice of the group which affords the given form.

 \begin{lemma}
   Let $h$ be a weight $k$  modular form for a finite index subgroup $\G$ of $SL_2(\Z)$ and $\G_h=\{\g\in \SL : h|\g=h\}$ which is a supergroup of $\G$. Then $h\in S_k^{\text{prim}}(\G)$ if and only if $h\in S_k^{\text{prim}}(\G_h)$.
 \end{lemma}
 \begin{proof}
   We first show that the set $\G_h\cdot \G^c=\{g_1g_2 : g_1\in \G_h, g_2\in \G^c\}=\G_h^c$. It is easy to see $\G_h\cdot \G^c\subset \G_h^c$. Let $H$ be any principal congruence subgroup contained in $\G^c$, then $\G_h\cdot H=\{g_1g_2 : g_1\in \G_h, g_2\in H\}$ is a group by the normality of $H$, and it is congruence and contained in $\G_h^c$. Thus $\G_h\cdot H=\G_h^c$ which implies $\G_h^c\subset \G_h\cdot \G^c$. It follows that $[\G_h^c:\G_h]=[\G^c:\G^c\cap \G_h]$. Thus any transversal of $ \G^c\cap \G_h$ in $\G^c$ is also a transversal of $\G_h$ in $\G_h^c$. So $\Tr_{\G_h}^{\G_h^c}h=0$ if and only if $\Tr_\G^{\G_c} h= \Tr_{\G^c \cap \G_h}^{\G_h^c}\Tr_{\G}^{\G^c \cap \G_h}h = [\G^c \cap \G_h : \G]\Tr_{\G_h \cap \G^c}^{\G^c}h = 0$ since $h$ is a cusp form for $\G^c \cap \G_h$.

 \end{proof}

  An immediate consequence is
\begin{cor} If $h\in S_k^{\text{prim}}(\G)$, then $h\in S_k^{\text{prim}}(G)$ for any finite index subgroup $G$ of $\G$.
\end{cor}

Next, assuming a form is primitive noncongruence, we derive a few primitive noncongruence conclusions on forms related to this form.

\begin{lemma}
  Let $\g\in \text{GL}_2(\Q)$ such that $\g^n\in \G$ for some positive integer $n$. If $h$ is primitive noncongruence, so is  $h|\g$.
\end{lemma}
\begin{proof}
  Let $G=\cap _{j=1}^n \g^{-j}\G\g^j$. It is normalized by $\g$ as $\g^n\in \G$. By Corollary 1, $h\in S_k^{\text{prim}}(G)$.  Further, $\g^{-1} G^c \g$ contains a congruence subgroup, hence $G^c \cap \g^{-1} G^c \g$ is a congruence subgroup containing $G$ and contained in $G^c$, thus it is equal to $G^c$. This means that $\g^{-1} G^c \g \supseteq G^c$, or equivalently, $ \g G^c \g^{-1} \subseteq G^c$. Thus $ \g G^c \g^{-1}$ is a congruence subgroup containing $\g G \g^{-1} = G$, hence is equal to $G^c$. We have shown that $G^c$ is also invariant under conjugation by $\g$.
Since $\g$ normalizes $G$,  $h|\g\in S_k(G)$.  We have $\Tr_G^{G^c} h|{\g} = (\Tr_G^{G^c} h)|{\g} = 0$.
\end{proof}

\begin{cor}\label{cor:1} Let $K$ be a positive integer.
   Suppose $h=\sum a(n)q^{n/\mu}$ is primitive noncongruence. Then the following two functions are also primitive noncongruence: (a) the subseries $\sum a(Kn)q^{Kn/\mu}$, and (b) $h\otimes\phi=\sum a(n)\phi(n)q^{n/\mu}$ for any Dirichlet character $\phi$ of conductor $K$.
\end{cor}
\begin{proof}
  Let $\g=\begin{pmatrix} 1 &  \mu/K\\ 0 & 1 \end{pmatrix}$. Note that $\g^{K}\in \G$. The claims follow from the above lemmas since    $$\sum a(Kn)q^{Kn/\mu}=\frac 1K\sum_{j=1}^K f|\gamma^j,$$  and as discussed in \cite{ALi78}, $$ h \otimes \phi = \frac{1}{{\bf g}(\phi^{-1})} ~\sum_{j \mod K, ~(j, K) = 1} ~\phi(j)^{-1}h| \g^j,$$
where ${\bf g}(\phi^{-1}) = \sum_{j \mod K, ~(j, K) = 1} ~\phi(j)^{-1} e^{2 \pi i j/K}$ is the Gauss sum attached to $\phi^{-1}$.
\end{proof}

\begin{lemma}\label{lem:3}
  If $h(z)$ is a primitive noncongruence cusp form of weight $k$, then so is $h(z/K)$ for any positive integer $K$.
\end{lemma}
\begin{proof} Let $\g = \begin{pmatrix} 1 & 0\\0 & K \end{pmatrix}$. Then $\g^{-1} \G_0(K) \g = \G^0(K)$,  where $\G_0(K)=\{\g\in \SL : \g\equiv \begin{pmatrix}
   *&*\\0&*
 \end{pmatrix} \mod K\}$ and $\G^0(K)$ is defined similarly using congruence with lower triangular matrices. Suppose that $h$ is a weight $k$ primitive noncongruence cusp form for $\G'$, then it is also primitive congruence for
$G = \G' \cap \G_0(K)$ by Corollary 1.  The form $h| \g(z) = K^{k/2}h(z/K)$ is a weight $k$ cusp form for $\g^{-1}G\g$, which is a finite index subgroup of $\G^0(K)$. The congruence closure $G^c$ of $G$ is contained in $\G_0(K)$ so that its conjugate $\g^{-1} G^c \g$ is a congruence subgroup contained in $\G^0(K)$. We claim that $\g^{-1} G^c \g$ is the congruence closure of $\g^{-1}G\g$. This is because
$$\g^{-1} G^c \g \supseteq (\g^{-1} G \g)^c \supseteq \g^{-1} G \g,$$ which is equivalent to
$$ G^c \supseteq \g(\g^{-1} G \g)^c \g^{-1} \supseteq G.$$ Since $\g(\g^{-1} G \g)^c \g^{-1}$ is congruence, it is equal to $G^c$ by the minimality of $G^c$. This proves the claim. Thus $\Tr_{\g^{-1}G\g}^{\g^{-1} G^c \g} h|{\g} = (\Tr_G^{G^c} h)|{\g} = 0$ if $ \Tr_G^{G^c} h = 0$.
\end{proof}

 Now we proceed to prove Theorem \ref{thm:main2}. Assuming $f \in S_k^{\prim}(\G)$, we shall derive a contradiction. Let $m$ be an integer such that $a(m)\neq 0$. By Corollary \ref{cor:1} and Lemma \ref{lem:3}, $h=\sum a(mn)q^{n/\mu}$ is a primitive noncongruence cusp form for a finite index subgroup  $\tilde{G}$. Let $K$ be the square of the product of all primes $\le P(m)$ in \S3, and let $\phi$ be a Dirichlet character of conductor $K$.
Then $h\otimes\phi$ is a primitive noncongruence form for the group $G = \cap _{j=1}^n \g^{-j}\tilde{G}\g^j$, where $\g=\begin{pmatrix} 1 &  \mu/K\\ 0 & 1 \end{pmatrix}$. It has Fourier expansion $(h \otimes \phi)(z) = \sum_{n \ge 1} a(nm) \phi(n) q^{n/\mu}\neq 0$. Moreover, by Lemmas \ref{lem:1} and \ref{lem:2},
it agrees with $a(m)(g \otimes \phi)(z/\mu)$, a cusp form for a congruence subgroup $\G'$. So it is
a cusp form for the intersection $G \cap \G'$, which is contained in the congruence subgroup $G^c \cap \G'$.
Note that $G \cap (G^c \cap \G') = G \cap \G'$ and there is a principal congruence subgroup $H'$ of $G^c \cap \G'$ such that $G \cdot H' = G^c$. %
 Thus $[G^c : G] = [G^c \cap \G' : G \cap \G']$ and consequently, $0 = \Tr_G^{G^c} (h \otimes \phi ) = \Tr_{G \cap \G'}^{G^c \cap \G'} (h \otimes \phi) = [G^c \cap \G' : G \cap \G'](h \otimes \phi)$ since $h \otimes \phi$ is a cusp form
for $\G'$ and hence for $G^c \cap \G'$. This then implies $h \otimes \phi = 0$, a contradiction.

\end{document}